\documentclass[12pt,reqno]{amsart}
\usepackage{amsmath,amssymb,amsfonts}
\usepackage{eucal,enumerate}

\numberwithin{equation}{section}

\newtheorem{lem}{Lemma}[section]
\newtheorem{cor}[lem]{Corollary}
\newtheorem{prop}[lem]{Proposition}
\newtheorem{thm}[lem]{Theorem}

\theoremstyle{remark}
\newtheorem{exam}[lem]{Example}
\newtheorem{prob}[lem]{Problem}

\newcommand\ts{\textstyle}

\renewcommand{\phi}{\varphi}                 
\renewcommand{\epsilon}{\varepsilon}
\newcommand\eset{\emptyset}                  
\renewcommand\emptyset{\varnothing}

\newcommand\inv{^{-1}}
\newcommand\transpose{^{\text{\rm T}}}
\newcommand\supp{\operatorname{supp}}
\newcommand\Lat{\operatorname{Lat}}
\newcommand\rk{\operatorname{rk}}
\newcommand\codim{\operatorname{codim}}
\newcommand\Nul{\operatorname{Nul}}
\newcommand\Row{\operatorname{Row}}
\newcommand\vol{\operatorname{vol}}

\newcommand\cH{\mathcal{H}}
\newcommand\cL{\mathcal{L}}
\newcommand\cM{\mathcal{M}}
\newcommand\cS{\mathcal{S}}

\newcommand\bbR{\mathbb{R}}

\newcommand\bbZ{\mathbb{Z}}

\newcommand\setm{\setminus}
\newcommand\0{\hat 0}
\newcommand\1{\hat 1}
\newcommand\Eta{\mathrm{H}}

\renewcommand\qedsymbol{\ensuremath{\blacksquare}}

\begin{document}

\title{The Number of Nowhere-Zero Flows 	\\
	on Graphs and Signed Graphs
}

\author{Matthias Beck}
\address{Department of Mathematics, San Francisco State University, San Francisco, CA 94132, U.S.A.}
\author{Thomas Zaslavsky}
\address{Department of Mathematical Sciences, Binghamton University, Binghamton, NY 13902-6000, U.S.A.}

\begin{abstract}
A nowhere-zero $k$-flow on a graph $\Gamma$ is a mapping from the edges of $\Gamma$ to the set $\{\pm1, \pm2, \dots, \pm(k-1)\} \subset \bbZ$ such that, in any fixed orientation of $\Gamma$, at each node the sum of the labels over the edges pointing towards the node equals the sum over the edges pointing away from the node.  
We show that the existence of an \emph{integral flow polynomial} that counts nowhere-zero $k$-flows on a graph, due to Kochol, is a consequence of a general theory of inside-out polytopes.  The same holds for flows on signed graphs.  We develop these theories, as well as the related counting theory of nowhere-zero flows on a signed graph with values in an abelian group of odd order.
Our results are of two kinds: polynomiality or quasipolynomiality of the flow counting functions, and reciprocity laws that interpret the evaluations of the flow polynomials at negative integers in terms of the combinatorics of the graph.
\end{abstract}



\thanks{Parts of this paper was written while the first author was a Robert Riley Assistant Professor at Binghamton University, SUNY. The author thanks the Department of Mathematical Sciences at Binghamton for its hospitality.}
\thanks{The research of the second author was partially supported by National Science Foundation grant DMS-0070729.}
\thanks{The authors are grateful for corrections and comments by Aaron Pixton and an anonymous referee.}

\date{Version of March 5, 2006. To appear in \emph{Journal of Combinatorial Theory Series B}}

\maketitle

\tableofcontents

\bigskip\hrule\bigskip

\emph{Note to publisher:}

This paper does NOT have a ``corresponding author'' or ``senior author''.

All authors are EQUAL.

All authors are able to answer correspondence from readers.

\bigskip

For editorial purposes ONLY, contact the writer of the cover letter of submission.

\bigskip\hrule\bigskip


\section{Introduction}  \label{flow}

A \emph{nowhere-zero flow} on a graph $\Gamma = (V,E)$, with values in an 
abelian group $A$, is a mapping $x: E \to A$ such that, for every node $v 
\in V$,
\begin{equation}\label{E:flow}
\sum_{h(e)=v} x(e) = \sum_{t(e)=v} x(e),
\end{equation}
and which never takes the value $0$. Here $h(e)$ and $t(e)$ are respectively the head and tail of the edge $e$ in a (fixed) orientation of $\Gamma$. (In a certain sense, described below, $x$ is independent of the chosen orientation.) A \emph{nowhere-zero $k$-flow} is an integral flow (i.e., $A=\bbZ$) whose absolute values are in $[k-1] := \{1,2,\ldots, k-1\}$.  
Nowhere-zero flows are nicely surveyed in \cite{Selected} and \cite{HandbookFlows}.

It has long been known that the number of nowhere-zero flows with values in a finite abelian group of order $k$ is a polynomial function of $k$ (Tutte \cite{Tring}).  Recently Kochol \cite{Kochol} discovered that the number of nowhere-zero $k$-flows is also a polynomial in $k$, although not the same polynomial.  Here we show that this fact is a consequence of a general theory of counting lattice points in inside-out polytopes \cite{IOP}.  
Furthermore, we extend Kochol's theorem in two ways: by a reciprocity law that combinatorially interprets negative arguments, and to signed graphs (in which each edge is positive or negative), where the polynomial becomes a quasipolynomial of period two: that is, a pair of polynomials, one for odd values of $k$ and the other for even $k$ (Theorem \ref{T:sgflowpoly}); 
and we partially extend to signed graphs Tutte's counting theory for nowhere-zero flows in abelian groups (Theorem \ref{T:sgmodflows}).  
Still further, the concept of reciprocity in lattice-point counting leads us to a geometrical interpretation of the number of totally cyclic orientations that are compatible with a given $k$-flow (Theorems \ref{T:flowpoly} and \ref{T:sgflowpoly} and Corollary \ref{C:sgtco}), a fact which parallels Stanley's theorem \cite{AOG} that the chromatic polynomial of $\Gamma$ evaluated at negative integers counts acyclic orientations compatible with a node-labelling of $\Gamma$.


\section{The method of inside-out polytopes}

The theory of inside-out polytopes was motivated by the problem of counting the integral points of a \emph{rational convex polytope} (the convex hull of finitely many rational points in $\bbR^d$) that do not lie in any of the members of a particular rational hyperplane arrangement.  
A (homogeneous, real) \emph{hyperplane arrangement} is a finite set of homogeneous hyperplanes in $\bbR^d$ (that is, hyperplanes that contain the origin); it is \emph{rational} if each of its hyperplanes is spanned by the rational points it contains or, alternatively, if each hyperplane has a rational normal vector.
Suppose we are given a rational convex polytope $P$ spanning $\bbR^d$ and a rational hyperplane arrangement $\cH$.  Then $(P,\cH)$ is a \emph{rational inside-out polytope}.  
More generally, $P$ and $\cH$ may lie in a rational (and homogeneous) subspace $Z$ that is spanned by $P$.  
We introduced inside-out polytopes recently \cite{IOP}; one motivation was to explain graph-theoretic enumeration functions geometrically.  We consider inside-out theory essential to understanding our results on integral flows. 
 
A \emph{region} (more precisely, an \emph{open region}) of $\cH$ is a connected component of $\bbR^d \setm \bigcup\cH$; its closure is a \emph{closed region}.  
The arrangement \emph{induced} by $\cH$ in a subspace $S$ of $\bbR^d$ is
$$
\cH^S := \{ H\cap S : H\in\cH, H\not\supseteq S \}.
$$
The \emph{intersection lattice} of $\cH$ is
$$
\cL(\cH) := \big\{ {\ts\bigcap} \cS : \cS \subseteq \cH \big\},
$$ 
ordered by reverse inclusion \cite{FUTA}; its elements are the \emph{flats} of $\cH$.  $\cL$ is a \emph{geometric lattice} with $\0 = \bigcap \eset = \bbR^d$ and $\1 = \bigcap \cH$.  (For matroids and geometric lattices we refer to \cite{Oxley} or \cite{EC1}.) 
The \emph{M\"obius function} of $\cL$ is the function $\mu: \cL \times \cL \to \bbZ$ defined recursively by 
\begin{equation*}
\mu(R,S) := \begin{cases}
	0			&\text{if } R \not\leq S, \\
	1			&\text{if } R = S, \\
	-\sum_{R \leq U < S} \mu(R,U)	&\text{if } R < S.
	\end{cases}
\end{equation*}
(Sources for the M\"obius function are, \emph{inter alia}, \cite{FCT} and \cite{EC1}.) 

\subsection{Ehrhart theory} 
We outline the classical Ehrhart theory of lattice-point enumeration in polytopes.  (See \cite{Ebook}, \cite[Section 4.6, pp.\ 235--241]{EC1}, or \cite{BR}.)  
We have a rational convex polytope $P$ that spans a subspace $Z$.  $P^\circ$ denotes the relative interior of $P$.   The \emph{denominator} of $P$ is the least common denominator of all the coordinates of vertices of $P$.  
We denote by $\vol{P}$ the \emph{volume} of $P$, normalized with respect to $Z \cap \bbZ^d$; that is, we take the volume of a fundamental domain of the integer lattice in $Z$ to be 1.  
To explain this last, we note that $Z\cap\bbZ^d$ is linearly equivalent to $\bbZ^{\dim Z} \subseteq \bbR^{\dim Z}$; a \emph{fundamental domain} is a domain in $Z$ that corresponds to the unit hypercube 
$[0,1]^{\dim Z} \subseteq \bbR^{\dim Z}$, under some invertible linear transformation that carries $Z\cap\bbZ^d$ to $\bbZ^{\dim Z}$.  
(See \cite[pp.\ 238--239]{EC1} or \cite[Section 5.4]{BR} for more detail.)  
When $Z = \bbR^d$ this is the ordinary volume.

A \emph{quasipolynomial} is a function $Q(t)=\sum_0^d c_i(t) \, t^i$ defined on $\bbZ$ with coefficients $c_i$ that are periodic functions of $t$.  Then $Q$ is a polynomial $Q_{\bar t}$ on each residue class $\bar t$ modulo some integer, called the \emph{period} of $Q$; these polynomials are the \emph{constituents} of $Q$. 
(For an introduction to quasipolynomials see, e.g., \cite[Chapter 4]{EC1} or \cite{BR}.)  

The subject of Ehrhart theory is the (closed) Ehrhart counting function
$$
E_P(t):= \# ( tP \cap \bbZ^d )
$$ 
and the open Ehrhart counting function, $E_{P^\circ}(t)$.  Ehrhart's theorem is that $E_P$ and $E_{P^\circ}$ are quasipolynomials with leading term $(\vol P) t^{\dim P}$ and with periods that divide the denominator of $P$.  It follows that one can define both counting functions for negative integers; the Ehrhart--Macdonald Reciprocity Theorem is that $E_{P^\circ}(t) = (-1)^{\dim P} E_P(-t)$.

\subsection{Inside-out Ehrhart theory} 
An \emph{open region of $(P,\cH)$} is a nonempty intersection with $P^\circ$ of an open region of $\cH$ (thus it is full-dimensional in the span of $P$).  A \emph{closed region of $(P,\cH)$} is the closure of an open region.  
A \emph{vertex} of $(P,\cH)$ is a vertex of any such region.  
The \emph{denominator} of $(P,\cH)$ is the least common denominator of all coordinates of all vertices. 

The fundamental counting functions associated with $(P,\cH)$ are two quasipolynomials: the (\emph{closed}) \emph{Ehrhart quasipolynomial},
\begin{align*}
E_{P,\cH}(t) &:= \sum_{x\in  t\inv \bbZ^d}m_{P,\cH}(x),
\intertext{where the \emph{multiplicity} $m_{P,\cH}(x)$ of $x\in \bbR^d$ with respect to $\cH$ and $P$ is defined through}
m_{P,\cH}(x) &:= \begin{cases}
\text{the number of closed regions of $(P,\cH)$ that contain $x$}, 
  & \text{if } x \in P, \\
0, & \text{if } x \notin P,
\end{cases}
\intertext{and the \emph{open Ehrhart quasipolynomial}, }
E_{(P,\cH)^\circ}(t) &:= \# \left( t\inv \bbZ^d \cap \left[ P^\circ \setminus \big( \bigcup\cH \big) \right] \right) .
\end{align*}
(This is simplified from the notation $E^\circ_{P^\circ,\cH}$ of \cite{IOP}.  The reader may think of $(P,\cH)^\circ$ as the set $P^\circ \setminus \big( \bigcup\cH \big)$, called the \emph{relative interior} of $(P,\cH)$, but that is not necessary in order to read this paper.)  
If $P$ spans a subspace $Z$, then $\cH$ also lies in $Z$ and the multiplicity of $x$ is defined to be 0 if $x \notin Z$.
The names of our counting functions are justified by the fact that in the absence of $\cH$ we recover classical Ehrhart theory, and by one of the main results in \cite{IOP}:

\begin{thm}[{\cite[Theorem 4.1]{IOP}}] \label{T:quasi}
If $(P,\cH)$ is a closed, full-dimensional, rational inside-out polytope in $Z \subseteq \bbR^d$, then $E_{P,\cH}(t)$ and $E_{(P,\cH)^\circ}(t)$ are quasipolynomials in $t$, with period equal to a divisor of the denominator of $(P,\cH)$, with leading term $(\vol{P}) t^{\dim P}$, and with constant term $E_{P,\cH}(0)$ equal to the number of regions of $(P,\cH)$.  Furthermore,
\begin{equation}\label{E:reciprocity}
E_{(P,\cH)^\circ}(t) = (-1)^{\dim P} E_{P,\cH}(-t).
\end{equation}
\end{thm}

In particular, if $(P,\cH)$ is integral then $E_{P,\cH}$ and $E_{(P,\cH)^\circ}$ are polynomials.  
The proof, though more general (and arrived at independently), is similar to Kochol's proof of Theorem \ref{T:flowpoly}(a): we sum the Ehrhart quasipolynomials of the pieces into which $\cH$ dissects $P$.

For the second theorem of \cite{IOP} we will use here, we need the notion of transversality: $\cH$ is called \emph{transverse} to $P$ if every flat $U \in \cL(\cH)$ that intersects $P$ also intersects $P^\circ$, and $P$ does not lie in any of the hyperplanes of $\cH$.

\begin{thm}[{\cite[Theorem 4.2]{IOP}}] \label{T:ehrhyp} 
If $P$ and $\cH$ are as in Theorem \ref{T:quasi}, then
\begin{equation}\label{E:oehrhyp}
E_{(P,\cH)^\circ}(t) = \sum_{U\in \cL(\cH)} \mu(\0,U) E_{P^\circ\cap U}(t),
\end{equation}
and if $\cH$ is transverse to $P$,
\begin{equation} \label{E:ehrhyp}
E_{P,\cH}(t) = \sum_{U\in \cL(\cH)}|\mu(\0,U)|E_{P\cap U}(t).
\end{equation}
\end{thm}

Transversality is always satisfied in the applications here.

\subsection{Matrix matroids} 
We shall want a general lemma about matroids of hyperplane arrangements induced by coordinate arrangements.  We start with the hyperplane arrangement $\cH_m$, consisting of the coordinate hyperplanes in $F^m$ for some field $F$, and we take any subspace $S$.  Then $S$ induces an arrangement $\cH_m^S$ in $S$.  
(We are interested in the reals, but the lemma is valid for any field.) 

Any homogeneous hyperplane arrangement $\cH$ has a matroid $M(\cH)$, whose ground set is the set of hyperplanes and whose rank function is $\rk \cS = \codim\left( \bigcap \cS \right)$ for $\cS \subseteq \cH$.  This matroid is simply the linear dependence matroid of the normal vectors of the hyperplanes.  
The \emph{column matroid} of a matrix $A$, $M(A)$, is the matroid of linear dependence of its columns; to keep the notation correct we take the ground set to be the set of indices of columns.  
The \emph{chain-group matroid} of a subspace $S \subseteq F^m$ is the matroid $N(S)$ on $[m]$ whose circuits are the minimal nonempty supports of vectors in $S$. 
$\Lat M$ denotes the lattice of closed sets of a matroid $M$.  Thus $\Lat M(\cH) \cong \cL(\cH)$.

We refer in the following lemma to orientations of oriented matroids.  The reader who is not familiar with oriented matroids should think of orientations of a graph or (later, in Section \ref{sg}) of a signed graph.

\begin{lem} \label{L:inducedhyp}
Let $A$ be an $n \times m$ matrix with entries in a field $F$, let $\cH_m = \{ H_e : e \in [m] \}$ be the arrangement of coordinate hyperplanes in $F^m$, and let $U = \Row A$, the row space, and $Z = \Nul A$, the null space.  
\begin{enumerate}

\item[(a)] The mapping $e \mapsto H_e \cap U$ is a matroid isomorphism from $M(A)$ to $M(\cH_m^U)$.  Also, $e \mapsto H_e \cap Z$ is a matroid isomorphism of $N(\Row A)$ with $M(\cH_m^Z)$.

\item[(b)] The mapping 
$$
S \mapsto E_U(S) := \{ e \in [m] : H_e \supseteq S \}
$$
is the isomorphism of $\cL(\cH_E^U)$ with $\Lat M(A)$ induced by the first mapping in (a).  The mapping
$$
S \mapsto E_Z(S) := \{ e \in [m] : H_e \supseteq S \}
$$
is the isomorphism of $\cL(\cH_E^Z)$ with $\Lat N(\Row A)$ induced by the second mapping in (a).

\item[(c)] If $F$ is an ordered field, then the regions of $\cH_m^U$ correspond bijectively to the acyclic orientations of the oriented matroid of columns of $A$, and those of $\cH_m^Z$ correspond to the totally cyclic orientations.

\end{enumerate}
\end{lem}

\begin{proof} 
(a) Each coordinate hyperplane $H_e$ equals $b_e^\perp$ for a basis vector $b_e$ along the $e$-axis.  
Therefore, $H_e \cap U$ is the orthogonal complement in $U$ of $a_e$, the orthogonal projection of $b_e$ into $U$.  
Let $B := \{ b_e : e \in [m] \}$.  We may take $A$ to be the matrix whose columns are the vectors $a_e$, since it has the same row space as the original matrix $A$.  
By linear duality, $M(A) \cong M(\cH_m^U)$ under the correspondence $e \mapsto H_e \cap U$.

We treat $\cH_m^Z$ by taking $A^*$, a matrix whose row space is $\Nul A$, and applying the first part to $\Row A^*$.  $N(\Row A)$ is dual to $M(A)$; hence it is $M(A^*)$ because $\Row A$ and $\Nul A$ are dual chain-groups.  
(See, e.g., Tutte \cite[Chapter VIII]{Tbook} on primitive chain-groups.)

(b) Obvious from (a) and the definitions.

(c) Write $\cM(A)$ for the oriented matroid.  The proof is as in part (a), but relying on the dual relationship between $\cM(A)$ and $\cM(A^*)$ \cite[Section 6.3(c)]{Ziegler}.  Our $\cM(A)$ is Ziegler's $\cM(\Nul A)$, so our $\cM(A^*) = \cM(\Row A)$, which is dual to $\cM(\Nul A)$ because the row and null spaces are orthogonal complements.
\end{proof}


\section{Flows on graphs}\label{graphs} 

A \emph{flow} on a graph $\Gamma$ with values in an abelian group $A$, called an \emph{$A$-flow}, is a function $x: E \to A$ which satisfies \eqref{E:flow} for every node $v \in V$ (so it is like a nowhere-zero flow, but the flow value zero is allowed).  
This definition requires that the edges be oriented in a fixed way. 
The orientation is arbitrary; it is an artifact of notation, and to overcome this artificiality we define, for an oriented edge $e$, $e\inv$ to be the same edge in the opposite orientation and $x(e\inv) := -x(e)$.  With this law for flows, the validity of Equation \ref{E:flow} is independent of the choice of the orientation of $\Gamma$.

Tutte \cite[Section 6]{Tring} proved that the number of nowhere-zero $A$-flows on $\Gamma$ is a polynomial in $|A|$, independent of the actual group.  
We shall write $\bar\phi_\Gamma$ for this polynomial and call it the \emph{(strict) modular flow polynomial} of $\Gamma$.  (Usually $\bar\phi_\Gamma$ is called just the ``flow polynomial'' but we need to distinguish it from other flow polynomials.)  
As Tutte showed, the modular flow polynomial $\bar\phi_\Gamma(k)$ is the evaluation $(-1)^{\xi(\Gamma)} t_\Gamma(0,1-k)$ of the Tutte polynomial of $\Gamma$.  ($\xi(\Gamma)$ is the cyclomatic number $|E| - |V| + c(\Gamma)$, where $c(\Gamma)$ is the number of connected components.)  
We should mention that the total number of $A$-flows is the simple polynomial $|A|^{\xi(\Gamma)}$.  
Tutte further proved (in \cite[pp.\ 83--84]{Tchromatic}, based on \cite[Theorem VI]{Temb}) that a nowhere-zero $\bbZ_k$-flow exists if and only if there is a nowhere-zero $k$-flow, a \emph{$k$-flow} being an integer-valued flow whose values all satisfy $|x(e)| < k$.  
However, the actual number of nowhere-zero $k$-flows for $k>0$, which we write as $\phi_\Gamma(k)$, does not equal the number of nowhere-zero $\bbZ_k$-flows and indeed was never known to be a polynomial until the recent work of Kochol \cite{Kochol}.  
Kochol employed standard Ehrhart theory combined with a special construction to prove this.  
We shall show that Kochol's theorem is a natural consequence of inside-out polytope theory and can be extended to a reciprocity theorem that interprets $\phi_\Gamma(k)$ at negative arguments. 

A \emph{circle} in a graph is a 2-regular connected subgraph.  
A \emph{cycle} (or \emph{directed cycle}) is a circle in which the edges are oriented in a consistent direction.  
An orientation of $\Gamma$ is \emph{acyclic} if it has no cycles and \emph{totally cyclic} if every edge lies in a cycle.  
We call a totally cyclic orientation $\tau$ and a flow $x$ \emph{compatible} if $x\geq 0$ when it is expressed in terms of $\tau$.  
Taking the standpoint of the flow $x$, the nonzero edge set $\supp x$ has a preferred orientation, the one in which $x \geq 0$ (we call this $\tau(x)$; note that it orients only $\supp x$) and the zero edges are free to take up any orientation that makes $\Gamma$ totally cyclic.  
An \emph{isthmus} of a graph is an edge whose deletion increases the number of connected components.  There is no totally cyclic orientation if $\Gamma$ has an isthmus.  

The \emph{real cycle space} $Z$ is defined in $\bbR^E$ by Equation \eqref{E:flow}.  
To this space $Z$ we associate the polytope and arrangement 
$$P := Z \cap [-1,1]^E, \qquad \cH := (\cH_E)^Z,$$
 where $\cH_E$ is the arrangement of coordinate hyperplanes in $\bbR^E$.  
A $(k+1)$-flow is then precisely a point $x \in Z \cap \bbZ^d$ such that $\frac{1}{k} x \in P$ and a nowhere-zero $k$-flow is just a point $x \in Z \cap \bbZ^d$ such that $\frac{1}{k} x \in P \setminus \left( \bigcup\cH \right)$.  Consequently,  
\begin{equation}\label{E:allflowdef}
\phi^0_\Gamma(k+1) = E_{P}(k)
\end{equation}
and
\begin{equation}\label{E:flowdef}
\phi_\Gamma(k) = E_{(P,\cH)^\circ}(k).
\end{equation}

\begin{thm}\label{T:flowpoly}  
Given: a graph $\Gamma$ and its real cycle space $Z$.
\begin{enumerate}

\item[(a)] (Kochol \cite{Kochol})  $\phi_\Gamma(k)$ is a polynomial function of $k$ for $k = 1,2,3,\ldots$\;.  It has leading term $(\vol{P}) k^{\xi(\Gamma)}$ if $\Gamma$ has no isthmi; otherwise it is identically zero.

\item[(b)]  Furthermore, $(-1)^{\xi(\Gamma)} \phi_\Gamma(-k)$ for $k\geq0$ equals the number of $(k+1)$-flows counted with multiplicity equal to the number of totally cyclic orientations of $\Gamma$ that are compatible with the flow.

\item[(c)] In particular, the constant term $\phi_\Gamma(0)$ equals the number of totally cyclic orientations of $\Gamma$, which equals $(-1)^{\xi(\Gamma)} \bar\phi_\Gamma(-1)$.

\item[(d)] Finally, the total number of $k$-flows for $k>0$, nowhere-zero or not, is a polynomial $\phi_\Gamma^0(k)$ satisfying $\phi_\Gamma^0(k) = (-1)^{\xi(\Gamma)} \phi_\Gamma^0(1-k)$, whose leading term is the same as that of $\phi_\Gamma(k)$ and whose constant term is $(-1)^{\xi(\Gamma)}$. 

\end{enumerate}
\end{thm}

We call $\phi_\Gamma$ the \emph{(strict) integral flow polynomial} of $\Gamma$ and $\phi_\Gamma^0(k)$ the \emph{weak integral flow polynomial}.  $\phi_\Gamma^0(k)$ already gives rise to a number of interesting computational problems, as discussed, for example, in \cite{Loera}. 

The reciprocity formula $\phi_\Gamma^0(k) = (-1)^{\xi(\Gamma)} \phi_\Gamma^0(1-k)$ leads to the conclusion that $\phi_\Gamma^0$ has no integral zeros.  For any $k>0$ there is an integral flow, the all-zero flow.  Reciprocity takes care of $k\leq0$.

\begin{proof}
For (a) we apply Theorem \ref{T:quasi} in $Z$.  
We call upon the total unimodularity of the matrix of the cycle equations \eqref{E:flow} to deduce that $P$ is a convex hull of integer lattice points.  Thus, $E_{(P,\cH)^\circ}$ is a polynomial, and as we saw in \eqref{E:flowdef}, it equals $\phi_\Gamma$.

Since $E_{P,\cH}(k)$ counts pairs $(x,R)$ where $x \in \bbZ^d \cap P$ and $R$ is a closed region of $\cH$ that contains $x$, part (b) follows if we show that the regions of $\cH$ correspond to the totally cyclic orientations of $\Gamma$ and a region of $\cH$ whose closure contains a chosen point $x \in Z \cap \bbZ^d$ corresponds to a totally cyclic orientation that is compatible with $x$.  The first statement was demonstrated by Greene and Zaslavsky in \cite[Section 8]{IWN}, based on the obvious bijection (given the fixed orientation of $\Gamma$) between orthants of $\bbR^E$ and orientations of $\Gamma$.  The second is then obvious.

Thus the constant term is the number of totally cyclic orientations.  The fact that this equals $t_\Gamma(0,2)$ is a theorem originally due to Las Vergnas (see \cite[Proposition 8.1]{CR} and \cite[remark after Theorem 1$'$, p.\ 296]{LasV}) and independently proved by Greene and Zaslavsky \cite[Corollary 8.2]{IWN}.

Part (d) is standard Ehrhart theory, because a $k$-flow is simply a point $x \in \bbZ^d$ such that $\frac{1}{k} x \in P^\circ$.  That is, $\phi_\Gamma^0(k) = E_{P^\circ}(k) = (-1)^{\xi(\Gamma)} E_{P}(-k)$ by Ehrhart reciprocity.  The constant term of $E_{P}(-k)$ is $1$, the Euler characteristic of $P$.  It is easy to see that $E_{P^\circ}(k) = E_{P}(k-1)$ for $k>0$.  Consequently, 
$$
\phi_\Gamma^0(k) = (-1)^{\xi(\Gamma)} E_{P}(-k) 
= (-1)^{\xi(\Gamma)} E_{P^\circ}(1-k) 
= (-1)^{\xi(\Gamma)} \phi_\Gamma^0(1-k) 
$$
if $k$ is a positive integer, whence for all $k$.
\end{proof}

\begin{prob} \label{Pr:flowvol}
Find a formula for, or a combinatorial interpretation of, the leading coefficient $\vol{P}$ of the integral flow polynomials.
\end{prob}

\begin{prob} 
Is there a combinatorial interpretation of $\bar\phi_\Gamma (-k)$ for 
$k \geq 2$? 
\end{prob}

\begin{exam}[Small graphs] \label{X:gflows}
We calculated the integral flow polynomials of some small graphs by counting integral $k$-flows on a computer for $1 \leq k \leq \xi(\Gamma)+2$ and interpolating to get the polynomial.  
The graphs were $mK_2$, the graph of $m$ parallel links, for $m = 3,4,5,6$, and $K_4$.  
We state our results along with the modular flow polynomials for comparison; the latter are 
$\bar\phi^0_\Gamma(k) = k^{\xi(\Gamma)}$ and 
$\bar\phi_\Gamma(k) = \chi_{\Gamma^*}(k)/k$, 
$\Gamma^*$ being the planar dual graph.  First, $3K_2$:
\begin{align*}
\bar\phi^0(k) &= k^2 , 	&\bar\phi(k) &= (k-1)(k-2) , 	\\
\phi^0(k) &= 3k^2 - 3k + 1 , 	&\phi(k) &= 3(k-1)(k-2) .
\end{align*}
Next, $4K_2$:
\begin{align*}
\bar\phi^0(k) &= k^3 , 	&\bar\phi(k) &= (k-1)(k^2-3k+3) , 	\\
\phi^0(k) &=  \frac{(2k-1)(8k^2-8k+3)}{3}, 
	&\phi(k) &= \frac{2(k-1)(8k^2-22k+21)}{3} .
\end{align*}
Next, $5K_2$:
\begin{align*}
\bar\phi^0(k) &= k^4 , \qquad	\bar\phi(k) = (k-1)(k^3-4k^2+6k-4) , 	\\
\phi^0(k) &= \frac{115k^4-230k^3+185k^2-70k+12}{12} , 	\\
\phi(k) &= \frac{5(k-1)(k-2)(23k^2-41k+36)}{12} .
\end{align*}
Next, $6K_2$:
\begin{align*}
\bar\phi^0(k) &= k^5 , \qquad \bar\phi(k) = (k-1)(k^4-5k^3+10k^2-10k+5) , \\
\phi^0(k) &=  \frac{2(2k-1)(44k^4-88k^3+71k^2-27k+5)}{10} , 	\\
\phi(k) &= \frac{(k-1)(176k^4-839k^3+1571k^2-1404k+620)}{10} .
\end{align*}
Finally, $K_4$:
\begin{align*}
\bar\phi^0(k) &= k^3 , 	&\bar\phi(k) &= (k-1)(k-2)(k-3) , 	\\
\phi^0(k) &= (2k-1)(2k^2-2k+1) , 	&\phi(k) &= 4(k-1)(k-2)(k-3) .
\end{align*}
\end{exam}

\begin{prob} \label{Pr:flowintegrality}
Is there any general reason why in some of these examples ($3K_2$ and $K_4$) both of the integral flow polynomials have integral coefficients and the integral nowhere-zero flow polynomial is a multiple of the modular polynomial?
\end{prob}

The totally cyclic orientations that are compatible with a flow $x$ are the totally cyclic extensions of $\tau(x)$ to $\Gamma$.  
We get such an extension by orienting the contraction $\Gamma/\supp x$ in a totally cyclic manner.  
The number of such orientations is given by Theorem \ref{T:flowpoly}(c).  Thus we have the following version of Theorem \ref{T:flowpoly}(b):

\begin{cor} \label{C:flowtc} 
$|\phi_\Gamma(-k)| = \sum_{x \in P \cap k\inv\bbZ^E}   |\bar\phi_{\Gamma/\supp x}(-1)|.$
\hfill\qedsymbol
\end{cor}

The \emph{multiplicity} of an integral point $x$ with respect to a rational hyperplane arrangement $\cH$ is the number of closed regions of $\cH$ that contain $x$.  
As we pointed out in proving Theorem \ref{T:flowpoly}(b), the number of closed regions of $\cH_E^Z$ that contain a flow $x$ is equal to the number of totally cyclic extensions to $\Gamma$ of $\tau(x)$.  Hence, the multiplicity of the integral flow $x$ with respect to $\cH_E^Z$ is another interpretation of the quantity $(-1)^{\xi(\Gamma/\supp x)} \bar\phi_{\Gamma/\supp x}(-1)$.

Other formulas arise from the intersection expansions of Theorem \ref{T:ehrhyp}, but as we need its M\"obius function, first we have to find the lattice $\cL(\cH)$ explicitly.  We do so in the more general context of signed graphs.


\section{Flows on signed graphs}  \label{sg} 

The best way to understand the cycle equations \eqref{E:flow} is in terms of the incidence matrix, which we expound in the general context of signed or bidirected graphs.  

Formally, a \emph{signed graph} $\Sigma = (\Gamma, \sigma)$ consists of a graph $\Gamma$ and a function $\sigma$ from the set of links and loops of $\Gamma$ to $\{+,-\}$.  (A \emph{link} has two distinct endpoints; a \emph{loop} has two coinciding endpoints.  In signed and bidirected graph theory it is convenient to have two more kinds of edges: a \emph{halfedge} has one endpoint and a \emph{loose edge} has no endpoints; neither of these has a sign.)  
If $T \subseteq E$, then $\Gamma|T$ or $\Sigma|T$ denotes the spanning subgraph whose edge set is $T$.
Each circle has a sign, which is the product of the signs of its edges.  A subgraph or edge set is called \emph{balanced} if it contains no halfedges and every circle in it has positive sign.  
(See \cite{NB} for the origin of signed graphs and balance; for the general theory of signed graphs see \cite{SG}.)

The \emph{bias matroid} (or \emph{signed-graphic matroid}) of $\Sigma$ \cite[Section 5 and Erratum]{SG}, written $G(\Sigma)$, can be defined by its rank function,
$$
r(T) = |V| - b(\Sigma|T) \text{ for an edge set } T,
$$
where $b(\Sigma|T)$ is the number of components of the subgraph $\Sigma|T$ that are balanced subgraphs, ignoring any loose edges.  
The circuits of $G(\Sigma)$ are of three kinds: a positive circle, a pair of negative circles that have a single common node, or a pair of node-disjoint negative circles together with a minimal connecting path; here one or both negative circles may be replaced by halfedges.  
A \emph{coloop} of $G(\Sigma)$ is an edge $e$ whose deletion makes an unbalanced component balanced; or which is an isthmus connecting two components of $\Sigma \setminus e$ of which at least one is balanced (all this by \cite[Theorem 5.l as corrected]{SG}).
We define the \emph{cyclomatic number} of $\Sigma$ to be $|E| - |V| + b(\Sigma)$.  
This is the rank of the dual $G^\perp(\Sigma)$ of the bias matroid.

A graph is \emph{bidirected} when each end of each edge is independently oriented.  
We express the bidirection by means of an \emph{incidence function} $\eta$ defined on the edge ends: the function is $+1$ if the arrow on that end points into the incident node, and $-1$ otherwise.  
For an edge end $\epsilon$, let $e(\epsilon)$ denote the edge and $v(\epsilon)$ the node incident to $\epsilon$.  We define 
$$
\eta(v,e) := \sum \big\{ \eta(\epsilon) : v(\epsilon) = v \text{ and } e(\epsilon) = e \big\}.
$$
Thus, $\eta(v,e) = 0$ if $v$ and $e$ are not incident or $e$ is a loose edge or positive loop.  
The \emph{bidirected incidence matrix} is
$$
\Eta(\Gamma,\eta) := \big( \eta(v,e) \big)_{V \times E}.
$$

A bidirection of a graph is really an orientation of a signed graph.  A link or loop $e$ with ends $\epsilon_1$ and $\epsilon_2$ has sign 
$$
\sigma(e) := - \eta(\epsilon_1)\eta(\epsilon_2).
$$  
In plain language, if the two arrows on $e$ point in the same direction, then $e$ is positive, but if they conflict, $e$ is negative.  We call $\eta$ an \emph{orientation of the signed graph $\Sigma$}.  
(See \cite{OSG}.  This notion corresponds to the ordinary notion of graph orientation if we identify an unsigned graph $\Gamma$ with the all-positive graph $+\Gamma$.)  
We call an \emph{(oriented) incidence matrix of $\Sigma$} any incidence matrix of an orientation of $\Sigma$, writing it $\Eta(\Sigma)$ (and $\Eta(\Gamma) := \Eta(+\Gamma)$).  
The ambiguity arising from the unspecified orientation can usually be ignored, since the only effect of changing the orientation is to negate some columns.  
Recall that reorienting $e$ to $e\inv$ also has the effect of replacing $x(e)$ by $x(e\inv) = -x(e)$ for every $x \in A^E$ ($A$ an abelian group).  
With these conventions we define a \emph{flow} on $\Sigma$ \emph{with values in $A$} as any $x \in A^E$ for which 
\begin{equation} \label{E:sgflow}
\Eta(\Sigma) x = 0,
\end{equation}
in other words, for which $x \in \Nul \Eta(\Sigma)$.  This definition generalizes that of a flow on a graph; naturally, then, we generalize theorems about flows.  
A \emph{$k$-flow} is an integral flow $x$ for which every $|x(e)| < k$, just as before.  
Nowhere-zero $k$-flows on signed graphs were studied by Bouchet \cite{Bouchet} and Khelladi \cite{Khelladi}, searching for the smallest possible $k$.

A \emph{cycle} in an oriented signed graph is a circuit such that no node has all arrows pointing into the node or all arrows pointing out of the node.  An orientation is \emph{acyclic} if it has no cycles and \emph{totally cyclic} if every edge belongs to a cycle.

Suppose we have a bidirected graph.  \emph{Switching} a node $v$ means changing $\eta$ to $\eta^v$ defined by
$$
\eta^v(\epsilon) = \begin{cases}
\eta(\epsilon)  &\text{ if } v(\epsilon) \neq v, \\
-\eta(\epsilon) &\text{ if } v(\epsilon) = v.
\end{cases}
$$
The associated switched signed graph is denoted $\Sigma^v$.  It is obtained from $\Sigma$ by negating all links incident with $v$.


\subsection{Group-valued flows}  \label{sgmodular} 

We begin our treatment of signed graphs with the analog of the modular flow polynomial, since as far as we know it has not been published.  
The \emph{Tutte polynomial} $t_\Sigma(x,y)$ of $\Sigma$ is defined to be that of $G(\Sigma)$.  (It is not equal to that of the underlying graph unless $\Sigma$ is balanced.  
See \cite{BOTutte} for the Tutte polynomial of a matroid.)

\begin{thm} \label{T:sgmodflows}
For each signed graph $\Sigma$ there is a polynomial $\bar\phi_\Sigma(k)$ such that the number of nowhere-zero flows on $\Sigma$ with values in a finite abelian group $A$ of odd order is $\bar\phi_\Sigma(|A|)$.  In fact, 
$$
\bar\phi_\Sigma(k) = (-1)^{\xi(\Sigma)} t_\Sigma(0,1-k).
$$
\end{thm}

\begin{proof}
For our proof by induction on $|E|$ we need deletion and contraction of edges in signed graphs, as in \cite{SG}.  
Deletion is simply removing the edge with no other change.  
The contraction $\Sigma/e$ by an edge $e$ is done differently for the four kinds of edge.  A loose edge or positive loop is simply deleted.  A halfedge or negative loop is deleted and its supporting node is also deleted, but any other edges incident to that node are retained; such an edge that was a link becomes a halfedge, but if it was a loop or halfedge it becomes a loose edge.  
A positive link is deleted and its endpoints are identified to a single node.  To contract a negative link we switch an endpoint so the link becomes positive; then the link is deleted and the endpoints are identified.

In the proof we fix $A$ and an orientation $\eta$ of $\Sigma$ and write $\bar\phi_\Sigma(A)$ for the number of nowhere-zero $A$-flows on $\Sigma$.

It is clear that if $\Sigma$ has components $\Sigma_1, \ldots, \Sigma_c$, then $\bar\phi_\Sigma(A) = 
\bar\phi_{\Sigma_1}(A) \cdots \bar\phi_{\Sigma_c}(A)$.  
Analogously, $t_\Sigma = t_{\Sigma_1} \cdots t_{\Sigma_c}$.  Therefore, we may assume $\Sigma$ is connected.

Suppose $e$ is a positive link of $\Sigma$, oriented from $v$ to $w$, such that deleting $e$ from $\Sigma$ does not increase the number of balanced components; that is, $e$ is not a coloop.  Let $x''$ be a nowhere-zero $A$-flow on $\Sigma/e$.  We can define a flow on $\Sigma$ by
\begin{align*}
x(f) &= x''(f) 	\qquad\text{ if } f \neq e, \\
x(e) &= \sum_{f\neq e} \eta(v,f) x(f).
\end{align*}
Then $\Eta(\Sigma)x=0$.  If $x(e)\neq 0$, $x$ is a nowhere-zero flow on $\Sigma$; if $x(e)=0$, $x$ is a nowhere-zero flow on $\Sigma \setminus e$.  Conversely, any nowhere-zero flow on $\Sigma \setminus e$ or $\Sigma$ gives rise to one on $\Sigma/e$.  
Since the assumption about $e$ guarantees that $\Sigma$, $\Sigma/e$, and $\Sigma \setminus e$ all have the same number of balanced components, 
\begin{equation}\label{E:flowdc}
(-1)^{\xi(\Sigma)} \bar\phi_\Sigma(A) = 
(-1)^{\xi(\Sigma/e)} \bar\phi_{\Sigma/e}(A) 
+ (-1)^{\xi(\Sigma \setminus e)} \bar\phi_{\Sigma \setminus e}(A).
\end{equation}

If $e$ is a negative link, subject to the same hypothesis on balance, with endpoints $v$ and $v'$, we switch $v$ to make $e$ positive.  Because a flow on $\Sigma$ remains a flow after switching, $\bar\phi_{\Sigma^v}(A) = \bar\phi_\Sigma(A).$  Thus, \eqref{E:flowdc} is valid for a negative link.

The analogous Tutte-polynomial formula, $t_\Sigma = t_{\Sigma/e} + t_{\Sigma \setminus e}$, is valid as long as $e$ is not a coloop in $G(\Sigma)$, but that is our balance assumption.  In this case, therefore, the theorem is valid for $\Sigma$ by induction.

If on the other hand an edge $e$ is a coloop, then $t_\Sigma(0,y)=0$ and, by the proofs of Bouchet's Lemmas 2.4 and 2.5 \cite{Bouchet}, also $\bar\phi_\Sigma(A) = 0$.  
(The oddness of $A$, so that $2a=0$ implies $a=0$ in $A$, is needed at this point.)

The preceding arguments reduce the theorem to the case of a signed graph with a single node $v$.  
The flow in a loose edge or positive loop can be any nonzero value in $A$, independent of all other flow values.  
Therefore, $\bar\phi_\Sigma(A) = (|A|-1)^l \bar\phi_{\Sigma_0}(A)$ where $l$ is the number of positive loops and loose edges and $\Sigma_0$ is $\Sigma$ with all such edges deleted.  
Now let us assume $\Sigma$ has one node and its edges are halfedges $e_1, \ldots, e_i$ and negative loops $f_1, \ldots, f_j$.  We may assume they are oriented into $v$, so the total inflow at $v$ is 
$$
x(e_1) + \cdots + x(e_i) + 2x(f_1) + \cdots + 2x(f_j) = 0
$$
by the cycle condition \eqref{E:sgflow}.  If $i+j=0$, there is one such flow, in agreement with $t_\Sigma = 1$.  
If $i+j=1$, there is none (again we require $|A|$ odd) so $\bar\phi_\Sigma = 0$, in agreement with $t_\Sigma(0,y)=0$.  
If $i+j>1$, we get a nowhere-zero $A$-flow by taking arbitrary nonzero values for $x$ except on one edge, say $x(e)$, which is determined by the rest through the cycle equations.  
The number of such choices in which $x(e)\neq0$ is $\bar\phi_\Sigma(A)$; the number in which $x(e)=0$ is $\bar\phi_{\Sigma\setminus e}(A)$.  Therefore,
\begin{align*}
\bar\phi_\Sigma(A) 
	&= (-1)^{i+j-1} (1-|A|)^{i+j-1} - \bar\phi_{\Sigma\setminus e}(A) \\
	&= (-1)^{\xi(\Sigma/e)}t_{\Sigma/e}(0,1-|A|) 
	 - (-1)^{\xi(\Sigma\setminus e)}t_{\Sigma\setminus e}(0,1-|A|)
\intertext{by the fact that $G(\Sigma/e)$ has rank 0 and induction on $i+j=|E|$,}
	&= (-1)^{\xi(\Sigma)}t_\Sigma(0,1-|A|)
\end{align*}
by the fact that $e$ is not a loop or coloop in $G(\Sigma)$.  Our theorem now follows by induction.
\end{proof}

\begin{prob} \label{Pr:sgfloweven}
Is there any significance to $\bar\phi^0_\Sigma$ or $\bar\phi_\Sigma$ evaluated at even natural numbers?
\end{prob}

Theorem \ref{T:sgmodflows} means there is a polynomial $\bar\phi_\Sigma(k)$, which we call the \emph{(strict) modular flow polynomial}, such that for any odd positive number $k$, $\bar\phi_\Sigma(k)$ is the number of nowhere-zero flows on $\Sigma$ with values in any fixed abelian group of order $k$.  
This is reminiscent of signed-graph coloring, where only odd values of $\lambda$ can be interpreted in the chromatic polynomial $\chi_\Sigma(\lambda)$ \cite{SGC}.  For coloring, though, there is another polynomial which counts restricted colorings when evaluated at even integers; and the two polynomials are related (as we showed in \cite{IOP}).  
We wonder whether there could be something similar with the modular flow polynomial, and whether flows and colorings might be connected through duality of signed graphs, analogously to the duality of colorings and flows on planar graphs.

\begin{exam}\label{X:sgfloweven}
The number of $A$-flows depends on the group $A$ when it has even order.  
Consider the signed graph consisting of two negative loops at a node $v$.  Orient both loops into $v$.  The cycle equation \eqref{E:flow} is $2x(e) + 2x(f) = 0$.  
We compare two groups of order 4.  
If $A$ is the Klein four-group $\bbZ_2 \times \bbZ_2$, then every $x : E \to A$ satisfies the cycle equations and we have 
$$
\bar\phi^0_\Sigma(A) = 16, \qquad \bar\phi_\Sigma(A) = 9.
$$
If $A = \bbZ_4$, then the $A$-flows are $(x(e),x(f)) \in \{0,2\}^2 \cup \{1,3\}^2$ and 
$$
\bar\phi^0_\Sigma(A) = 8, \qquad \bar\phi_\Sigma(A) = 5.
$$
For comparison, at odd values of $k$ we have 
$$
\bar\phi^0_\Sigma(k) = k, \qquad \bar\phi_\Sigma(k) = k-1.
$$
\end{exam}

\begin{cor} \label{C:sgtco}
The number of totally cyclic orientations of $\Sigma$ equals $(-1)^{\xi(\Sigma)} \bar\phi_\Sigma(-1)$.
\end{cor}

\begin{proof}
The number of totally cyclic reorientations of an orientation of a matroid $M$ is $t_M(0,2)$ \cite{CR,LasV}.  Since cycles in an orientation of $\Sigma$ are the same as cycles in the corresponding orientation of $G(\Sigma)$ \cite{OSG}, the number of totally cyclic orientations of $\Sigma$ equals $t_{G(\Sigma)}(0,2) = (-1)^{\xi(\Sigma)} \bar\phi_\Sigma(-1)$.  
\end{proof}


\subsection{Integral $k$-flows on signed graphs}  \label{sgintegral} 

Now it is time for integral flows.  For $k>0$ let 
$$
\phi_\Sigma(k) := \text{ the number of nowhere-zero $k$-flows on }\Sigma.  
$$
As with abelian-group flows, $\phi_\Sigma = 0$ if there is a coloop in $G(\Sigma)$.  Let 
$$
\phi_\Sigma^0(k) := \text{ the number of all $k$-flows on }\Sigma , 
$$
also for $k>0$.  We take $Z$ to be the real cycle space $\Nul \Eta(\Sigma)$---this is the solution space of Equation \eqref{E:flow}---and, just as with unsigned graphs, 
$$P := Z \cap [-1,1]^E \quad\text{ and }\quad \cH := \cH_E^Z$$
 where $\cH_E$ is the coordinate-hyperplane arrangement.  
 As with ordinary graphs, a flow $x$ and an orientation $\eta$ are \emph{compatible} if $x\geq0$ when expressed in terms of $\eta$.

\begin{thm} \label{T:sgflowpoly}
\begin{enumerate}

\item[(a)] 
For any signed graph $\Sigma$, $\phi_\Sigma(k)$ is a quasipolynomial function of $k$ for $k = 1,2,3,\ldots$\;.  Its period is 1 or 2, and is 1 if $\Sigma$ is balanced.  $\phi_\Sigma(k)$ has leading term 
$(\vol{P}) k^{\xi(\Sigma)}$ if $G(\Sigma)$ has no coloops; otherwise $\phi_\Sigma(k)$ is identically zero.

\item[(b)] Furthermore, $(-1)^{\xi(\Sigma)} \phi_\Sigma(-k)$ for $k\geq0$ equals the number of $(k+1)$-flows counted with multiplicity equal to the number of compatible totally cyclic orientations of $\Sigma$.

\item[(c)] In particular, the constant term $\phi_\Sigma(0)$ equals the number of totally cyclic orientations of $\Sigma$, which equals $(-1)^{\xi(\Sigma)} \bar\phi_\Sigma(-1)$.

\item[(d)] Finally, $\phi_\Sigma^0(k)$ is a quasipolynomial of period 1 or 2 (period 1 if $\Sigma$ is balanced) whose leading term is the same as that of $\phi_\Sigma(k)$ and whose constant term is $(-1)^{\xi(\Gamma)}$.  Furthermore, $\phi_\Sigma^0(k) = (-1)^{\xi(\Gamma)} \phi_\Sigma^0(1-k)$.

\end{enumerate}
\end{thm} 

From (d) we see that, even when $\phi^0$ has period $2$, in a way it consists only of one polynomial.  The odd constituent, $\phi^0_{\text{odd}}$, is determined by the even constituent through
$$
\phi^0_{\text{odd}}(2j+1) = (-1)^{\xi(\Gamma)} \phi^0_{\text{even}}(-2j) .
$$

\begin{prob}\label{Pr:period}
Is there a nontrivial unbalanced signed graph for which $\phi_\Sigma$ has period $1$?  (For the dual function, the chromatic quasipolynomial, there is not.  See the second remark after \cite[Corollary 5.9]{IOP}.)  The graph of one node in Example \ref{X:sgflows} shows that trivial examples exist.  As for $\phi^0_\Sigma$, there exists at least one such graph with more than one node, but still very small; see Example \ref{X:sgflows}.

The mere existence of such examples is significant: it suggests that the duality between flows and colorings is imperfect and limited.
\end{prob}

\begin{lem} \label{L:cyclehalfint}
The vertices of $(P,\cH)$ are half integral.
\end{lem}

\begin{proof}
We may as well assume $\Sigma$ has a halfedge at every node; thus $\Eta(\Sigma)$ contains an identity matrix $I_n$.  A vertex is a solution of $\Eta(\Sigma)x = 0$ with $|E| - n$ coordinates of $x$ set equal to fixed values in $\{0, 1, -1 \}$.  Let $B$ be the edge set whose coordinates in $x$ are left undetermined, let $B^c := E \setminus B$, and write $x = (x_B, x_{B^c})\transpose$.  Then $x$ is the unique solution of $\Eta(\Sigma|B)x_B = - \Eta(\Sigma|B^c)x_{B^c}$.

The remainder of the proof is based on work of Jon Lee.  Lee proved that the null space $\Nul \Eta(\Sigma)$ is 2-regular \cite[Proposition 9.1]{Leesub} and that if $A$ is a nonsingular square matrix for which $\Nul [I \mid A]$ is 2-regular, then $A\inv b$ is half integral for every integral vector $b$ (a special case of \cite[Proposition 6.1]{Leesub}; the definition of 2-regularity need not concern us).  These facts applied to $A = \Eta(\Sigma|B)$ imply that the solution of $\Eta(\Sigma|B) x = b$ is half integral for any $b \in \bbZ^B$.  Apply this to $b = - \Eta(\Sigma|B^c)x_{B^c}$.
\end{proof}

\begin{proof}[Proof of Theorem \ref{T:sgflowpoly}]
The proof is similar to that of Theorem \ref{T:flowpoly}.  In (a) and (d), instead of total unimodularity we have Lemma \ref{L:cyclehalfint} to tell us that the denominator of $(P,\cH)$, hence the period of the Ehrhart quasipolynomials, divides 2.

For (b) we need to show that the regions of $\cH$ correspond to the totally cyclic orientations of $\Sigma$.  The latter are the totally cyclic reorientations of the natural orientation of $G(\Sigma)$, which is the oriented matroid of columns of $\Eta(\Sigma)$ \cite[Theorem 3.3]{OSG}.  Now we apply Lemma \ref{L:inducedhyp}(b).

For (c) we use Corollary \ref{C:sgtco}.
\end{proof}

\begin{exam}\label{X:sgflows}
We calculated the flow polynomials and modular flow polynomials of some small signed graphs.  
First we treat the signed graph with two negative loops at one node, for comparison with the modular flow polynomials in Example \ref{X:sgfloweven}:
\begin{equation*}
\phi^0_{\Sigma}(k) = 2k-1 , \qquad \phi_{\Sigma}(k) = 2(k-1) .	
\end{equation*}
Another signed graph, only slightly larger, has two nodes joined by positive and negative edges and at each node either a halfedge or a negative loop.  We write $\pm K_2^{(i,j)}$ for this graph if $i$ nodes have halfedges and $j$ other nodes have negative loops and $\phi_{(i,j)}$ for the polynomials; the examples we calculated are where $i+j=2$.  The modular nowhere-zero flow polynomial is the characteristic polynomial of the dual matroid of $G(\pm K_2^{(i,j)})$, which is the four-point line.
The modular polynomials are
\begin{align*}
\bar\phi^0_{{(i,j)}}(k) &= k^2 ,  \\
\bar\phi_{{(i,j)}}(k) &= (k-1)(k-3) .
\end{align*}
The integral polynomials are
\begin{align*}
\phi^0_{{(0,2)}}(k) = \phi^0_{{(2,0)}}(k) &= 2k^2-2k+1 , \\	
\phi_{{(0,2)}}(k) = \phi_{{(2,0)}}(k) &= \begin{cases}
2(k-1)(k-3)	&\text{if $k$ is odd},	\\
2(k-2)^2	&\text{if $k$ is even}	\end{cases} 
\end{align*}
and
\begin{align*}
\phi^0_{{(1,1)}}(k) &=  \begin{cases}
\frac{1}{2}(3k^2-2k+1)	&\text{if $k$ is odd},	\\
\frac{1}{2}(3k^2-4k+2)	&\text{if $k$ is even},	\end{cases}	\\
\phi_{{(1,1)}}(k) &= \begin{cases}
\frac{1}{2}(3k^2-12k+9)	&\text{if $k$ is odd},	\\
\frac{1}{2}(3k^2-14k+16)	&\text{if $k$ is even}.	\end{cases}
\end{align*}
\end{exam}

\begin{prob} \label{Pr:sgflows}
We do not understand why $\pm K_2^{(0,2)}$ and $\pm K_2^{(2,0)}$ have the same flow quasipolynomials but $\pm K_2^{(1,1)}$ does not, nor why $\phi^0_{{(0,2)}}$ has period one, nor why the odd constituent of $\phi_{{(0,2)}}$ is an integral multiple of $\bar\phi_{{(i,j)}}$.
\end{prob}

\subsection{Half integrality and the incidence matrix} \label{half}

In our geometric treatment of graph coloring in \cite[Section 5]{IOP} we noticed an important half-integrality property of the signed-graphic hyperplane arrangement $\cH[\Sigma]$, which consists of the dual hyperplanes to all the columns in the incidence matrix of $\Sigma$.  
It is the next lemma, which involves a polytope $P=[0,1]^d$ and the affine arrangement $\cH''[\Sigma]$ which is $\cH[\Sigma]$ translated by the vector $\frac12(1,1,\ldots,1)\transpose$.

\begin{lem}[{\cite[Lemma 5.7]{IOP}}]  \label{L:sgdenominator}  
If $\Sigma$ is a signed graph, $(P,\cH''[\Sigma])$ has half-integral vertices.
\end{lem}

Despite their appearance, not only are Lemmas \ref{L:cyclehalfint} and \ref{L:sgdenominator} similar in statement, both also involve the incidence matrix.  
A \emph{basic solution} of a linear program $\Eta(\Sigma)\transpose x \leq b$, $x \geq 0$ is the solution of $\Eta(\Sigma|B)\transpose x = b$ for some set $B$ of $n$ linearly independent columns of $\Eta(\Sigma)$.  

\begin{prop} \label{L:sghalfint}
Every basic solution of $\Eta(\Sigma)\transpose x \leq b$, $x \geq 0$ is half integral for any $b \in \bbZ^n$.
\end{prop}

\begin{proof} [Proof that Lemma \ref{L:sgdenominator} and Proposition \ref{L:sghalfint} are equivalent.]
We may assume that $\Sigma$ contains a negative loop at every node.  
Write $a_e$ for the column of the incidence matrix that belongs to $e$.  
The equations of the hyperplanes in $\cH''[\Sigma]$ are 
$a_e\transpose x = \frac12 a_e\transpose (1,1,\ldots,1)\transpose$; 
this is 0 or 1 if $e$ is a loop or link but it is $\frac 12$ if $e$ is a halfedge at node $v_i$ and in that case the equation simplifies to $x_i=\frac12$.  
In the situation of Lemma \ref{L:sgdenominator}, the vertices of $(P,\cH''[\Sigma])$ are obtained by choosing from among the equations $x_i=0$, $x_i=1$ (from $P$) and $x_i=\frac12$, $a_e\transpose x = 0$ or $1$ (from $\cH''[\Sigma]$) a subset of $n$ equations whose coefficient matrix is nonsingular.  
An equation $x_i = 0$ or 1 has the form $a_e\transpose x = 0$ or 1 with $e$ the negative loop at $v_i$, so we need not consider it separately.  Thus, half integrality of vertices of $([0,1]^n, \cH''[\Sigma])$ implies half integrality of the solution of $\Eta(\Sigma|B)\transpose x = b$ for any $b \in \bbZ^n$.  
The converse is obvious.
\end{proof}

Our attention was drawn to the relationship between Lemma \ref{L:sgdenominator} and half integrality (and Lee's work) by a recent manuscript of Appa and Kotnyek \cite{Appa}.


\subsection{Nowhere-zero flows reduce to flows with M\"obius complications} \label{mobius} 

The last main result expresses the nowhere-zero integral flow polynomial in terms of the weak integral flow polynomials of subgraphs.  We begin with structural lemmas.  As before, $Z$ is the cycle space, $P := [-1,1]^E \cap Z$, and $\cH := \cH_E^Z$.  For a flat $S \in \cL(\cH)$ we define 
$$
E(S) := \{ e \in E : S \subseteq H_e \} = \{e \in E: x(e) = 0 \text{ for all } x \in S \}.
$$
This is the $E_Z(S)$ of Lemma \ref{L:inducedhyp}(b).  We see that $E(S)^c$ is the union of the suppports of the vectors in $S$.

\begin{lem} \label{L:cyclelattice}
The lattice of flats of $\cH$ is isomorphic to the lattice of closed sets of the dual of the bias matroid, $G^\perp(\Sigma)$.  The isomorphism is given by $S \mapsto E(S)$.  The corresponding matroid isomorphism $G^\perp(\Sigma) \cong M(\cH)$ is given by $e \mapsto H_e \cap Z$.
\end{lem}

\begin{proof} 
This is an application of Lemma \ref{L:inducedhyp}.  The matrix is $\Eta(\Sigma)$.  We know $M(\Eta(\Sigma))$ is the bias matroid $G(\Sigma)$ by \cite[Theorem 8A.1]{SG}, so $G^\perp(\Sigma)$ is the 
chain-group matroid of $\Row \Eta(\Sigma)$.  The lemma applies because the real cycle space $Z = \Nul \Eta(\Sigma)$.
\end{proof} 

(Lemma \ref{L:cyclelattice} holds good for the canonical hyperplane representation of any $F^*$-gain graph $\Phi$, for any field $F$ \cite[Section 2]{BG4}.  Denoting an incidence matrix by $\Eta(\Phi)$ and taking $Z = \Nul \Eta(\Phi)$, we have $M(\cH_E^Z) = G^\perp(\Phi)$.  But we digress.)

\begin{lem} \label{C:noflow}
$\phi_\Sigma$ is identically zero if and only if $G(\Sigma)$ has a coloop.
\end{lem}

\begin{proof} 
This is an application of Bouchet's theorem on integral chain-group matroids \cite[Proposition 3.1]{Bouchet}: the chain-group has a nowhere-zero chain if and only if the matroid has no coloop.  In our case the chain-group is the group of integral flows, $\Nul \Eta(\Sigma)$.  
Its chain-group matroid is dual to that of $\Row \Eta(\Sigma)$, which is dual to the column matroid of $\Eta(\Sigma)$, which is $G(\Sigma)$.
\end{proof} 

\begin{lem} \label{L:subflow}
A flat $S$ of $\cH$ can be represented as $[\Nul \Eta(\Sigma|E(S)^c)] \times \{0\}^{E(S)}$.
\end{lem}

\begin{proof}
The lemma is obvious from the definitions of $Z$ and $E(S)$.
\end{proof}

\begin{thm} \label{T:flowhyp}
Take a signed graph $\Sigma$.  Letting $T$ range over all subsets of $E$, or merely over all for which $G(\Sigma)|T$ has no coloops (that is, all complements of flats of the dual bias matroid $G^\perp(\Sigma)$),
\begin{equation} \label{E:flowhyp}
\phi_\Sigma(-k) = \sum_T |\mu(\0,T^c)| \phi^0_{\Sigma|T}(k+1)
\end{equation}
and
\begin{equation} \label{E:nzflowhyp}
\phi_\Sigma(k) = \sum_T \mu(\0,T^c) \phi^0_{\Sigma|T}(k+1),
\end{equation}
where $\mu$ is the M\"obius function of $G^\perp(\Sigma)$ and 
$\0$ is the set of coloops of $G(\Sigma)$.
\end{thm}

\begin{proof} 
The polytope and the arrangement are transverse because $\bigcap \cH$ intersects $P^\circ$.  

Since $\phi_{\Sigma|T} = 0$ if $G(\Sigma)|T$ has a coloop by Lemma \ref{C:noflow}, the two ranges of summation are equivalent.  For $S \in \cL(\cH)$, by Lemma \ref{L:subflow} we know that
$$
P \cap S = [-1,1]^E \cap \big( [\Nul \Eta(\Sigma|E(S)^c)] \times \{0\}^{E(S)} \big).
$$
Take $T = E(S)^c$.  Then
$$
P \cap S = \big( [-1,1]^T \cap Z' \big) \times \{0\}^{T^c}
$$
where $Z'$ is the real cycle space of $\Sigma|T$.  Its Ehrhart polynomial equals $\phi^0_{\Sigma|T}(k+1)$.  

Now the result follows from Lemma \ref{L:cyclelattice}, Theorem \ref{T:ehrhyp}, Equation \eqref{E:flowdef}, and Theorem \ref{T:flowpoly}(d).
\end{proof} 

A reminder: to apply the theorem to a graph $\Gamma$, take $\Sigma = +\Gamma$.

Thus, the strict integral flow polynomial can be expressed in terms of the weak polynomial and invariants of $G^\perp(\Sigma)$.  If the weak polynomial were as simple as in the case of colorings, where it is a 
monomial \cite[Section 5]{IOP}, we would have a nice formula for the number of nowhere-zero $k$-flows.  But such is not the case.

It may be helpful to list some characterizations of the edge sets that support nowhere-zero integral flows.

\begin{prop} \label{L:cyclesubspaces}
For $T \subseteq E := E(\Sigma)$, the following properties are equivalent.
\begin{enumerate}
\item[(i)] $G(\Sigma)|T$ has no coloops.
\item[(ii)] $\Sigma|T$ has a totally cyclic orientation.
\item[(iii)] $\Sigma|T$ has a nowhere-zero integral flow.
\item[(iv)] $\Sigma|T$ has a nowhere-zero real flow.
\item[(v)] $T^c$ is closed in the dual bias matroid $G^\perp(\Sigma)$.
\item[(vi)] $T = E(S)^c$ for some flat $S \in \cL((\cH_E)^Z)$.
\end{enumerate}
\end{prop}

\begin{proof} 
(i) $\iff$ (ii) for graphs is Robbins' theorem \cite{Robbins}.  We could prove it for signed graphs, but the simplest approach is via oriented matroids.  We know that the number of totally cyclic reorientations of an orientation of a matroid $M$ is $t_M(0,2)$ \cite{CR,LasV} and that this equals 0 if and only if $M$ has a coloop.  Apply that to the natural orientation of $G(\Sigma)$.

(i) $\iff$ (iii) by Lemma \ref{C:noflow}.

(iii) $\implies$ (iv) is trivial.

(iv) $\implies$ (i) by the proofs of \cite[Lemmas 2.4 and 2.5]{Bouchet}, which amount to saying that any flow on $\Sigma$ with values in an abelian group where $2a = 0 \implies a = 0$ must be zero on every coloop.  Here the group is the additive group of $\bbR$.

(v) $\iff$ (i):  By matroid duality the complements of the closed sets in $G^\perp(\Sigma)$ are the edge sets that do not contain a coloop of $G(\Sigma)$.

(v) $\iff$ (vi):  This is Lemma \ref{L:cyclelattice}.
\end{proof}


\end{document}